\documentclass[10pt]{article}

\usepackage{fullpage}

\usepackage{paralist}
\usepackage{wrapfig,subfigure}
\usepackage{amssymb,amsmath,amstext,amsgen,amsthm,amsbsy,amsopn,amsfonts,mathrsfs,graphicx,framed}
\usepackage{setspace}
\usepackage{tikz}
\usetikzlibrary{shapes,snakes,shadows,calc,positioning,backgrounds}
\definecolor{dgreen}{rgb}{0,.7,0}

\usepackage{hyperref}
\hypersetup{colorlinks,breaklinks,linkcolor=black,urlcolor=black,anchorcolor=black,citecolor=black}

\theoremstyle{definition}

\usepackage{bm}

\usepackage{natbib,listings}
\usepackage{epstopdf}
\usepackage{algorithm}

\newcommand{\e}{\ensuremath{\varepsilon}}

\title{Perturbation Theory\\ \normalsize Entry for the \emph{Enclyclopedia of Theory} with Sage Publications.\\[.5cm] Current word count: 4416 (Target 4000)}

\author{Nicolas Fillion \\ Dept. of Philosophy\\ SFU \and Robert M. Corless \\Dept. of Applied Mathematics\\UWO}

\begin{document}

\maketitle

%\tableofcontents

%\pagebreak

This article aims to explain essential elements of perturbation theory and their conceptual underpinnings. It is not meant as a summary of popular perturbation methods, though some illustrative examples are given to underline the main methodological insights and concerns. We also give brief explications of the mathematical notions of \emph{limit}, \emph{continuity}, \emph{differentiability}, \emph{convergence}, and \emph{divergence}, which provide the necessary foundation.

\section{Introduction}

The role of perturbation theory is to study the effects of small disturbances (or perturbations) on the behaviour of a (possibly only hypothetical) system, irrespective of the source of the disturbance.  For instance, the disturbance may have its physical origin in the environment of a system, as when a system receives small amounts of energy from the outside, or when an experimental setup is slightly disturbed by the vibration caused by a truck on a nearby street. The disturbance may have been added to the modelling equations by design in order to simplify the study of a system, as is the case when we neglect small quantities to obtain more mathematically tractable equations. The disturbance may even be due to experimental error and uncertainty or to computational error.

The effects of small disturbances are sometimes small too, in which case they are usually called \emph{regular} perturbations. As \citet[p.~1]{SimmondsMann} explain, regular perturbations ``are assumed nearly every time we construct a mathematical model of a real world phenomenon. Our choice of language reflects this: the flow is \emph{almost} steady, the density varies \emph{essentially} with altitude only, the conductivity is \emph{virtually} independent of temperature, the spring is \emph{nearly} linear, the friction is \emph{practically} negligible.'' However, mother nature isn't always that cooperative and sometimes the effect of small disturbances is very large indeed; in such cases we talk of \emph{singular} perturbations. In addition to being essential to the mathematical disciplines, this distinction has become an important topic of discussion in philosophy of science; see, e.g., \citet{Batterman(2002)}.\nocite{FillionMoir(2018)} We will return to the distinction between regular and singular perturbations later. %and discuss more precise definitions.

Perturbation theory is central to the scientific method in two important but distinct respects. On the one hand, good applied mathematical practice demands that one considers possible disturbances before assessing whether a mathematical model correctly captures the behaviour of a system. Thus, even ``exact'' mathematical theories and models enter the realm of perturbation theory as soon as they are applied. On the other hand, while it is typical of introductory courses on ``mathematical $X$''---e.g., mathematical physics, biology, economics, etc.---to focus exclusively on problems that can be solved exactly, this desirable event of complete tractability only rarely happens in applied mathematical practice, either because we don't know exactly solvable models or because the accurate models turn out to be theoretically impossible to solve exactly. Hence instead of using ``demonstrative solution methods'' one necessarily relies on methods that extract approximate solutions. In this case, a problem is said to be `solved' in the sense that we can get approximations that are close enough, even if we can never actually solve it exactly (we return to this below).

One of the beautiful and powerful aspects of perturbation theory is that it addresses the two concerns with the very same set of mathematical methods. These methods, to be sure, differ from other mathematical methods that are central to science. For instance, deductive or axiomatic approaches focus on whether axioms are true and on whether certain statements are logical consequences of sets of axioms. Probability theory concentrates on quantifying how probable it is that a statement will be true. However, the few remarks above suggest that in most concrete contexts, it is not enough to ask what is true, what logically follows from what, and what is probable; we also need to ask: \emph{if things were changed or disturbed a little, what would the effect be}? This is precisely the role of perturbation theory and what distinguishes it from other mathematical disciplines.

\section{Mathematical background}

To understand the main concepts and methods of perturbation theory, we introduce a few essential mathematical ideas. We believe that they will be accessible to most readers who have an intuitive grasp of the concept of derivative and, to a lesser extent, of the concept of integral.

We begin with the crucial notion of \emph{limit} and the related notion of \emph{convergence}. Intuitively, a limit is the mathematical object (often, a number or a function) that is approached by an infinite \emph{sequence}. Consider a classical example. Even if we know by definition that $\pi$ is the ratio of the circumference of a circle to its diameter, it is impossible to write down its value as a number in conventional format since its decimal expansion is infinite and non-repeating. However, around 250BC, Archimedes made the first known use of sequences to calculate the value of $\pi$ to arbitrary precision, at least in principle.
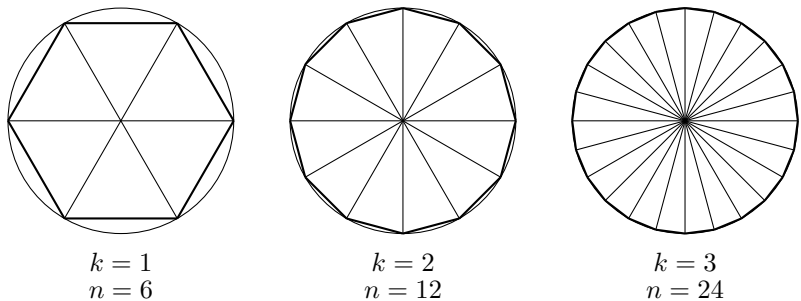
\begin{figure}
\centering
\begin{tikzpicture}[scale=1.5]
\begin{scope}[shift={(0,0)}]
  \draw (0,0) circle (1cm);
\def\n{6};
\foreach \k in {1,...,\n}	
	{
	\draw ({cos(\k*360/\n)},{sin(\k*360/\n}) edge[thick] ({cos((\k+1)*360/\n)},{sin((\k+1)*360/\n});
	\draw (0,0) -- ({cos(\k*360/\n)},{sin(\k*360/\n});
	}
\draw (0,-1.25) node {$k=1$};
\draw (0,-1.5) node {$n=\n$};
\end{scope}

\begin{scope}[shift={(2.5,0)}]
  \draw (0,0) circle (1cm);
\def\n{12};
\foreach \k in {1,...,\n}	
	{
	\draw ({cos(\k*360/\n)},{sin(\k*360/\n}) edge[thick] ({cos((\k+1)*360/\n)},{sin((\k+1)*360/\n});
	\draw (0,0) -- ({cos(\k*360/\n)},{sin(\k*360/\n});
	}
\draw (0,-1.25) node {$k=2$};
\draw (0,-1.5) node {$n=\n$};
\end{scope}

\begin{scope}[shift={(5,0)}]
  \draw (0,0) circle (1cm);
\def\n{24};
\foreach \k in {1,...,\n}	
	{
	\draw ({cos(\k*360/\n)},{sin(\k*360/\n}) edge[thick] ({cos((\k+1)*360/\n)},{sin((\k+1)*360/\n});
	\draw (0,0) -- ({cos(\k*360/\n)},{sin(\k*360/\n});
	}
\draw (0,-1.25) node {$k=3$};
\draw (0,-1.5) node {$n=\n$};
\end{scope}
% Archimedes had shown 223/71 < � < 22/7, which is for n=
\end{tikzpicture}
\caption{The first three items in Archimedes' sequence of polygons circumscribed in the unit circle to approximate $\pi$.}
\label{archimedes}
\end{figure}
As we see in figure \ref{archimedes}, the perimeter of the  hexagon circumscribed in the unit circle is equal to 6, since it can be decomposed in six equal equilateral triangles of side 1. Archimedes' brilliant idea was that, if we know the side $s$ of a polygon circumscribed in a circle, we can find the side $s'$ of a polygon with twice the number of sides by elementary trigonometry (see figure \ref{trig}):
\begin{figure}
\centering
\begin{tikzpicture}[scale=5]
\draw (1,0) arc (0:50:1);
\draw (0,0) -- node[below] {1} ({cos(0)},{sin(0)}) -- ({cos(25)},{sin(25)}) -- node[pos=0.03,below] {$\alpha$} node[midway,below] {$\beta$}  (0,0) -- ({cos(50)},{sin(50)}) -- node[near start,below left] {$\frac{s}{2}$} node[near end, below left] {$\frac{s}{2}$} (1,0);
\draw[very thick] ({cos(50)},{sin(50)}) -- ({cos(25)},{sin(25)});
\draw (1.4,.5) node {$\beta = \sqrt{1-\left(\frac{s}{2}\right)^2}$};
\draw (1.4,.3) node {$\alpha= 1-\beta$};
\draw (1.4,.1) node {${\bf s'} = \sqrt{\alpha^2+\left(\frac{s}{2}\right)^2}$};
\end{tikzpicture}
\caption{The trigonometric relation exploited by Archimedes.}
\label{trig}
\end{figure}
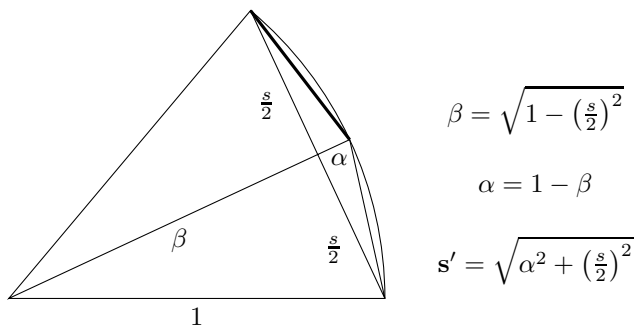
\begin{align}
s' = \sqrt{\left(1-\sqrt{1-\left(\frac{s}{2}\right)^2}\right)^2+\left(\frac{s}{2}\right)^2}
\end{align}
Using this formula, %we find that the successive approximations of $\pi$ are
%\[
%p_k = \frac{1}{2}(6\cdot2^{k-1})s_k = \frac{1}{2}k\sqrt{\left(1-\sqrt{1-\left(\frac{s_{k-1}}{2}\right)^2}\right)^2+\left(\frac{s_{k-1}}{2}\right)^2}
%\]
%with $s_1=1$.
we can successively find the length of the sides of circumscribed polygons with $6,12,24,48,\ldots$ sides, from which it is easy to find approximations of the ratio of the perimeter to the diameter. This gives us the following sequence of approximations:
\[
\begin{array}{c@{\quad}c@{\quad}c@{\quad}c}
p_1=3.0000 & p_2=3.1058 & p_3=3.1326 & p_4=3.1394 \\ p_5=3.1410 & p_6=3.1415 & p_7=3.1416 & \ldots
\end{array}
\]
Here, we round to the fourth decimal place for convenience. The bigger $k$ is, the closer $p_k$ will be to $\pi$. In that sense, we say that the sequence of $p_k$ converges to $\pi$ as $k$ goes to infinity (denoted $k\to\infty$). Alternatively, we say that $\pi$ is the limit of this sequence as $k\to \infty$. This example illustrates convergence. Archimedes in fact demonstrated more, namely the existence of the unique number that we call $\pi$, by demonstrating that the perimeters of inscribed polygons ended to the same number. We could have used this example to illustrate perturbation theory itself, by noting that each $p_k$ is a constructible perimeter of a perturbed figure.

%s(1)=1;
%p(1)=6;
%sp = @(s) sqrt( (1- sqrt(1-(s/2)^2) )^2 + (s/2)^2  );
%for k=2:5
%s(k)=sp(s(k-1));
%p(k)=3*2^(k-1)*s(k);
%end

Limits of sequences are also essential to the manipulation and interpretation of \emph{infinite series}:
\begin{align}
S = \sum_{k=0}^\infty a_n = a_0+a_1+a_2+\cdots + a_n + \cdots
\end{align}
Consider the perfectly well-defined finite sum
\begin{align}
S_N = \sum_{k=0}^N a_k = a_0 + a_1 + a_2+\cdots+a_N
\end{align}
of $N+1$ terms. The sums $S_1,S_2,S_3,\ldots$ form a sequence. The sum $S$ of an infinite series is interpreted as the limit of this sequence, i.e., the number or function that is approached as the number $N+1$ of terms summed approaches infinity. In mathematical notation, we write
\begin{align}
S = \lim_{N\to\infty} S_N = \lim_{N\to\infty} \sum_{k=0}^N a_k = \sum_{k=0}^\infty a_k\>.
\end{align}
If no such finite $S$ exists, the series is said to \emph{diverge}. However, if the series has a limit $S$, then it is possible to use a finite number $N$ of terms to approximate the value of $S$. Since this practice involves retaining the first $N$ terms of the series and discarding the rest, this is called a \emph{truncated series approximation}.

There is a substantive arsenal of tests available to determine whether series converge. But what is the status of series that don't converge to any number or function? In the nineteenth century, Abel provocatively put it as follows:
\begin{quote}
The divergent series are the invention of the devil, and it is a shame to base on them any demonstration whatsoever. By using them, one may draw any conclusion he pleases and that is why these series have produced so many fallacies and paradoxes [\ldots]. \citep[cited in][p.~218]{Hoffman(1998)}
\end{quote}
From the point of view of analysis, they are said to have no meaning since they don't refer to anything; alternatively the series is said not to be valid.

However, as we will see, there are important differences between this perspective and that which provides the conceptual foundations for perturbation theory. %In the next section, we will see explain how it is that series can be tremendously useful whether or not they are convergent.
Indeed, additional subtleties come into play when we consider \emph{series of functions} that may (or may not) converge to a function. For sequences of functions one has to be concerned with the convergence being \emph{uniform}; we will keep discussion of uniform convergence to a minimum but note its importance both in theory and in practice. This being said, many of the most successful applications of mathematics in science have only become possible by an improved understanding of how divergent series can be used to obtain reliable results.

First, we observe that functions themselves may have a limit near a point; this is an application of the same idea we've introduced above. Consider a function $f(x)$ defined on an interval (or more generally, on an open set) $A$, and suppose we wish to know what value the function $f(x)$ approaches as the argument $x$ becomes increasingly closer to a point $a$ within the interval. We can think of a sequence of arguments $a_0,a_1,a_2,a_3,\ldots$ and of the corresponding sequence of values $f_1=f(a_1),f_2=f(a_2),f_3=f(a_3),\ldots$. If the sequence of $f_k$s approaches a given finite value $L$ as the $a_k$s approach $a$, for any choice of sequence of $a_k$s that approach $a$, then we say that $L$ is the limit of the function $f(x)$ as $x$ approaches $a$. In rigorous mathematical notation, we say that a function $f(x)$ has the limit $L$ as $x\to a$ if
\begin{align}
(\forall \e>0)(\exists\delta>0)(x\in A\ \&\ |x-a|<\delta \Rightarrow |f(x)-L|<\e)\>.
\end{align}
A function $f(x)$ is said to be \emph{continuous} at the point $x=a$ provided that
\begin{align}
\lim_{x\to a} f(x) = f(a)\>.
\end{align}
It is continuous on an interval (or open set) $A$ if it is continuous at all points $a$ in the interval. This definition captures the intuitive idea that a function is continuous if its graph is an unbroken curve, without holes or jumps.

The treatment of limits at $\infty$ is nearly the same, because the ideas are very similar (some would say identical except in form). We say that
\begin{align}
\lim_{x\to\infty} f(x)= L
\end{align}
if
\begin{align}
 (\forall \e>0)(\exists X)(x>X \Rightarrow |f(x)-L|<\e)\>.
\end{align}
Here the large $X$ plays the same role as the small $\delta$ earlier.

This is known as the ``epsilon-delta'' definition of limit and it is due to the efforts of Cauchy, Bolzano, and Weierstrass to provide rigorous foundations for calculus in the early nineteenth century. This definitional approach provides concrete meanings for more intuitive phrases, used since the time of Newton, such as ``we can make the error $|f(x)-L|$ as small as we please by taking $x$ sufficiently close to $a$,'' without appealing to infinities and infinitesimals. See \citet{Boyer(1949)} for a nice history of calculus and \citet{Bell(2005)b} for a history of infinitesimals. It must be noted that many new approaches have revived the use of infinitesimals, such as non-standard analysis \citep[see, e.g.][]{Robinson(1966)} and smooth infinitesimal analysis \citep[see][]{Bell(2008)}.

Following Carath\'eodory, we say that $f$ is \emph{differentiable} at $x=a$ if there exists a function $\varphi(x)$ that is continuous at $x=a$ such that
\begin{align}
f(x) = f(a) + (x-a)\varphi(x)\>,
\end{align}
in which case $f'(a) = \lim_{x\to a} \varphi(x) = \varphi(a)$. This is of course the equation of the tangent of the function $f$ at the point $x=a$. The process can be repeated to find the derivative of the derivative (known as the second derivative) $f''(x)$, the third derivative $f'''(x)$, and so on. This gives rise to an \emph{expansion} of the function $f$ about the point $a$. When a function has $N+1$ derivatives in a neighbourhood of $x=a$, then there exists a continuous function $\varphi_{N+1}(x)$ such that
\begin{align}
f(x) &%= f(a) + f'(a)(x-a) + \frac{f''(a)}{2}(x-a)^2 +\frac{f'''(a)}{6}(x-a)^3 + \cdots \nonumber \\ &
= \sum_{k=0}^N \frac{1}{k!}f^{(k)}(a) (x-a)^k + \varphi_{N+1}(x) (x-a)^{N+1}\>. \label{taylor}
\end{align}
This is known as Taylor's theorem, and the resulting Taylor series are among the most useful in applied mathematics.

To simplify expressions such as the one in equation \eqref{taylor}, we can use the ``big-Oh'' notation, first used by Bachmann and popularized by Landau. To indicate that $f(x)$ and $g(x)$ are of similar size near $x=a$, we say that $f$ is order $g$ (denoted $f(x)=O(g(x))$ near $a$. More precisely, we say that
\begin{align}
f(x)=O(g(x)) \qquad \textrm{near } x=a
\end{align}
if there exists a constant $K$ such that for all $x$ sufficiently near $a$ we have
\begin{align}
|f(x)| \leq K|g(x)|.
\end{align}
In theory $K$ may be very large---indeed arbitrarily large, or arbitrarily small. But in practice, when the use of a perturbation method is successful, $K$ is usually of quite moderate size. % In an engineering sense, $f(x)=O(g(x))$ means that $f$ and $g$ are about the same size; the mathematical sense is as previously stated and the two agree only if $\frac{1}{10}\lesssim K \lesssim 10$, or so.
With this notation, equation \eqref{taylor} simply becomes
\begin{align}
f(x) & = \sum_{k=0}^N \frac{1}{k!}f^{(k)}(a) (x-a)^k + O\left( (x-a)^{N+1}\right)\>.
\end{align}

Following the same idea, but generalizing it, perturbation theory uses truncated series of functions:
\begin{align}
S_N(\e) = \sum_{n=0}^N a_n\phi_n(\e)\>. \label{perturbseries}
\end{align}
The sequence $\phi_0(\e),\phi_1(\e),\ldots,\phi_N(\e),\phi_{N+1}(\e),\ldots$ may be infinite; these functions are often called \emph{gauge functions}.
The most usual gauge functions are simple powers of $\e$ (i.e., $1,\e,\e^2,\ldots$) when we are interested by the behaviour of the function as $\e\to 0$, or powers of $(\e-a)$ when we are interested in the behaviour of the function as $\e\to a$. %\footnote{Expansion at infinity is scarcely more difficult. One uses gauge functions that vanish ever more rapidly at infinity, and one typically uses another letter (say $t$).}
 Notice that when $\e$ is near 0 or $a$, as the case may be, the higher the power $N$ of $\e^N$ or $(\e-a)^N$ is, the closer it will be to $0$; in other words, for gauge functions, each function in the sequence is smaller than the previous one near the limit point:
\begin{align}
\lim_{\e\to0} \frac{\phi_{n+1}(\e)}{\phi_n(\e)} = 0
\end{align}
This is only useful if the $\phi_k$ are somehow simpler than the original function $f$. Indeed we know an example,
\begin{align}
\mathrm{ln} z = \sum_{k=1}^N W^{\langle k\rangle}(z) + O(W^{\langle N+1\rangle}(z))\qquad \textrm{as } z\to\infty
\end{align}
of an asymptotic expansion of a simple function in terms of iterates of a less well-known function, the Lambert $W$ function, but if this is useful at all the example would be a startling exception \citep[see][]{Corless(1996)}. The normal case has each $\phi_k$ simpler than the original $f$.

One might also say that perturbation methods are an application of the theory of truncated \emph{formal} series. %In a typical application a finite number of terms of the series are used (indeed usually just two terms, or three). The purpose is to give a simple quantitative explanation of the effects of varying a parameter (often called $\varepsilon$ nowadays as the dominant mathematical culture since Cauchy uses $\varepsilon$ to denote a small, variously specified, number.)
 %
 %Because typically only a finite number of terms are used in perturbation expansions, say $N+1$ terms as in
% \begin{align}
%S_N(\e) = \sum_{n=0}^N a_n\phi_n(\e)\>,
%\end{align}
%the series are said to be ``truncated''.
 If one does not explicitly truncate the series as in equation \eqref{perturbseries} but writes $a_0\phi_0(\e)+a_1\phi_1(\e)+a_2\phi_2(\e)+\cdots$, %or even $\sum_{k\geq0} a_k\phi_k(\e)$
  one describes the series as ``formal''. This means that issues of convergence as $N\to\infty$ are usually irrelevant, and in practice divergent series are often used. The point is, and it is a simple one although there is a great deal of confusion in the literature, that  there are two limits of interests with series:
\[
\lim_{N\to\infty} S_N(\e) \>, \quad \textrm{where here \e\ is fixed, and } \e\neq0\>,
\]
encountered first in most calculus classes, whence issues of convergence of a sum of an infinite number of terms arise, and the other
\[
\lim_{\e\to0} S_N(\e)\>, \textrm{where here $N$ is fixed},
\]
the relevant limit for perturbation series. In the quite common case that the second limit as $\e\to0$ exists and is finite but the first does not (the series diverges as $N\to\infty$) the series is said to be `asymptotic'. The cure for divergence as $N\to\infty$ is simply to keep $N$ fixed. The other limit, as $\e\to0$, often provides all the accuracy one needs.

This point of view allows one to use a fixed number of terms (\emph{any} fixed number of terms) of a series that would diverge if we were so foolish as to let $N\to\infty$. %To fix notation suppose we have $N+2$ gauge functions $\phi_0(\e),\phi_1(\e),\phi_2(\e),\ldots,\phi_N(\e),\phi_{N+1}(\e)$ or if the limit $t\to\infty$ is used, $\{\phi_k(t)\}_{k=0}^{N+1}$. Suppose also that
%\begin{align}
%\lim_{\e\to0} \frac{\phi_{k+1}(\e)}{\phi_k(\e)}=0\quad k=0,1,\ldots,N
%\end{align}
%or the equivalent as $t\to\infty$ for $\phi_k(t)$. That is, each
Using gauge functions, i.e., functions that are smaller than the previous ones in a neighbourhood of the limit point, we find that the perturbation series for $f(a+\e)$ is
\begin{align}
f(a+\e) = \sum_{k=0}^N p_k\phi_k(\e) + O(\phi_{N+1}(\e))
\end{align}
if
\begin{align}
\lim_{\e\to0} \frac{|f(a+\e) - \sum_{k=0}^m p_k\phi_k(\e)|}{|\phi_{m+1}(\e)|}
\end{align}
is bounded (\textit{mutatis mutandis} as $t\to\infty$), for each $m=0,1,2,\ldots,N$. Again, we stress that $N$ is fixed and does not go to infinity.

\section{The use of series in perturbation theory}

With the mathematical concepts introduced in the last section, we are now ready to examine how perturbation methods work. When facing a hard problem, one can use the following conceptually (and with some practice, mathematically) simple three-step method.%\footnote{For a similar pedagogically insightful presentation, see \citet{Bender}.}

The \emph{first step} consists in introducing a small positive parameter $\e$, known as a perturbation parameter, in the hard problem that you may not know how to solve, or that can be provably unsolvable by explicit methods. So, instead of having only one hard problem, we are now in the presence of an infinite collection of hard problems, one for each value of \e. This may not look like progress, but that is only the first step. The idea of introducing a perturbation $\varepsilon$ is sometimes an obvious thing to do (e.g., when a value of $\varepsilon$ such as $0$ really simplifies things), and sometimes it is not obvious at all. Very clever choices are sometimes needed, and knowing what choice to make is in many respects an art rather than a science. However, in scientific use, the modelling context may sometimes dictate where to introduce the \e\ as it will correspond to a quantity that is small, such as air resistance, friction, etc. %But one important factor that comes into the choice of where to introduce $\varepsilon$ is that the \emph{un}perturbed problem, with $\varepsilon=0$, can be solved exactly.

For the \emph{second step}, assume that the answer is a function of $\varepsilon$ and assume that the answer is in the form of a truncated power series in $\varepsilon$:
\begin{align}
A\!N\!S(\varepsilon) = \sum_{n=0}^N a_n\e^n
\end{align}
More generally, we could use any collection of gauge functions instead of the powers of \e. This practice, which seems to assume that we have a valid series expansion without actual  mathematical evidence, is sometimes objected to; but we will explain below the sense in which it is legitimate. The main problem consists in calculating the coefficients $a_n$, and our examples below show how this can be achieved recursively. This is why perturbation theory is sometimes defined as follows:
\begin{quote}
Perturbation theory is a large collection of iterative methods for obtaining approximate solutions to problems involving a small parameter $\varepsilon$. \citep[p.~319]{Bender(1978)}
\end{quote}
Note also that since we aim to obtain a truncated series approximation, we will fix the number of coefficients that we will solve for ($N$ will often be only 1 or 2, in which case we will say respectively that we obtained a first-order or second-order solution).

In the third and final step, we set $\varepsilon$ to the value that recovers the original problem we wanted to solve (often, it will be $\e=1$) and just calculate the approximate solution by summing a finite number of terms of the perturbation series. That simple method very often works, but as we will see sometimes more care is needed \citep[for a more technically precise discussion and more examples, see the tutorial in][]{CorlessFillion(2019)}. We begin with an example in which it turns out to work well.

Consider the equation $x^5+x-1=0$. Our problem is to find a real root of this equation, i.e., a value of $x$ that satisfies this equation. The problem is hard, but we need not despair as we can use the method outlined above. Firstly, we introduce an \e\ as follows:
\begin{align}
x^5+\e x-1=0\label{algeq1}
\end{align}
We will choose to keep only two terms to obtain a first-order approximation. We thus assume that our truncated series solution will have the form $x(\e) = a_0 +a_1\e$. We substitute this expression in equation \eqref{algeq1}:
\begin{align}
(a_0 +a_1\e)^5 +\e (a_0 +a_1\e) - 1 =0
\end{align}
We expand the fifth-power term using Pascal's triangle:
\begin{align}
(a+b)^5 = a^5+5a^4b+10a^3b^2+10a^2b^3+5ab^4+b^5
\end{align}
After gathering like terms, we obtain
\begin{align}
a_0^5-1+ (5a_0^4a_1+a_0)\e + (10a_0^3a_1^2+a_1)\e^2 + O(\e^3) = 0\>, \label{matchpowers}
\end{align}
hiding the complicated higher-order terms in the $O$-symbol (we do not need them anyway, since the $\e^2$ term will already fail to be $0$ with this first-order approximation).

Because asymptotic power series are unique, an asymptotic power series equal to zero must have each coefficients equal to 0 \citep[see, e.g.][p.~6]{Jeffreys(1962)}. This fact justifies equations \eqref{cond1} and \eqref{cond2} below, and convergence as $N\to\infty$ is not needed. Thus, we then match the coefficients of matching powers of \e\ on the left-hand side and the right-hand side.  Since all coefficients of \e\ on the right-hand side are $0$, we obtain the following  equations for the first two coefficients:
\begin{align}
a_0^5 - 1 &= 0\label{cond1} \\
5a_0^4a_1+a_0 &=0 \label{cond2} %\\
%10a_0^3a_1^2+a_1 &=0 \label{cond3}
\end{align}
Since we were interested with the real root of this equation, we obtain the value $a_0=1$ from equation \eqref{cond1} and $a_1=-1/5$ from equation \eqref{cond2}. We see that the solution of the coefficients is a recursive process that is in itself a purely mechanical task that presents no essential mathematical difficulty, provided the zeroth order equation, which is generally non-linear, can be solved.

Accordingly, our approximate solution, which we denote $z(\e)$ in contrast with the exact but unknown solution $x(\e)$, can be written as follows:
\begin{align}
z(\e) = 1-\frac{\e}{5}
\end{align}
Now, if we set $\e=1$ in equation \eqref{algeq1}, we obtain our original unperturbed problem. Thus, our approximate root will be $z(1) = 1-1/5 = 0.8$; it does not compare too badly with the exact solution, which is $0.755$ rounded to three digits. If we were to do a bit more work, we would find that $a_2=-1/25$; accordingly our second-order approximation would be $0.76$, even better!

Notice that we have provided no evidence to justify any step in which we assumed that the root $x(\e)$ of equation \eqref{algeq1} can be meaningfully written as an infinite power series. But in fact, we have made no such assumption! The approach here can be justified by an explicitly finite approach borrowed from numerical analysis that in many cases allows physical interpretation, based on the concept of  \emph{residual} \citep[see, e.g.,][]{CorlessFillion(2013)}.

The idea is quite general; the unknown function $f(\e)$ usually satisfies a known equation $\Phi(f;\e)=0$ by definition (above, we had $\Phi(f;\e)=\Phi(x,\e)=x^5+\e x-1=0$). However, if we evaluate $\Phi$ using $z$ instead of $x$, we wouldn't obtain $0$ since $z$ isn't the exact solution. The residual $\Delta$ is the quantity we obtain once we substitute the exact function $f(\e)$ by its truncated series approximation:
\begin{align}
\Delta := \Phi\left(\sum_{k=0}^N a_k\e^k,\e\right)
\end{align}
Then $z$ satisfies
\[
\Phi(z;\e)-\Delta=0\>,
\]
not $\Phi(z;\e)=0$. As we see, the approximate solution exactly satisfies a modified equation.

In general if we correctly identify the first $N+1$ coefficients we obtain a residual that satisfies
\begin{align}
|\Delta|=O(\e^{N+1}) \label{resid}
\end{align}
as $\e\to0$. As we have explained, the symbol $O(\e^{N})$ is conveniently used to refer to a function that decays like $\e^N$ without having to specify it explicitly.  %This technique is equivalent to knowledge of the forward error if one knows the ``condition number'', namely the effect of small changes of the input to $\Phi$ on its solution. This is of course itself a first-order perturbation question, so we have the interesting spectacle of a theory recursively applied to itself!
In the example above, we implicitly used the notion of residual to find coefficients that guarantee that equation \eqref{resid} will be satisfied, i.e., that the residual will be of higher order than the terms that are part of our approximate solution. Instead of writing equation \eqref{matchpowers}, we could more exactly have written our equation for the residual,
\begin{align}
\Delta = a_0^5-1+ (5a_0^4a_1+a_0)\e + (10a_0^3a_1^2+a_1)\e^2 + O(\e^3)\>,
\end{align}
and by an iterative inspection of the formula we find that choosing $a_0=1$ guarantees that the residual is $O(\e)$ irrespective of the other coefficients, and that in addition choosing $a_1=-1/5$ guarantees that the residual is $O(\e^2)$. %\marginpar{I need to say somewhere that as $e\to0$, larger powers of $e^k$ are smaller!}
For small \e, this $z$ is the exact solution of nearly the right equation. %Because this happens to be a well-conditioned problem, this approximation is also nearly the right root: we have $|z-x|=O(\e^2)$ as well.

By doing more work, but using no new ideas, we may find $z_N = a_0+a_1\e+a_2\e^2+\cdots+a_N\e^N$ for any fixed $N\geq 1$ with $\Delta_N = z_N^5+\e z_N-1=O(\e^{N+1})$, which can be made arbitrarily small by taking $\e$ arbitrarily close to 0. Also, for this example, more is true. We may in this case take $N\to\infty$, and find an explicit formula for all $a_k$, thereby giving an infinite series---which turns out to be convergent---and the exact root. But that is no longer perturbation theory.

%The discussion above generalizes in a number of important ways. Firstly, perturbation theory is not restricted to power series (expansion in terms of powers of $\e$). Indeed any convenient collection of functions $\phi_0(\e),\phi_1(\e),\phi_2(\e),\ldots$ can be used, so long as each is smaller than the previous near the limit point:
%\begin{align}
%\lim_{\e\to0} \frac{\phi_{n+1}(\e)}{\phi_n(\e)} = 0
%\end{align}
%Functions with this property are known as \emph{gauge functions}. Perturbation theory is also not restricted to expansion about $\e=0$: any finite point can be considered simply by translation: $f(x)$ near $x=a$ is the same as $f(a+\e)$ near $0$. Expansion at infinity is scarcely more difficult. One uses gauge functions that vanish ever more rapidly at infinity, and one typically uses another letter (say $t$). An instance is
%\begin{align}
%y(t) = \frac{-1}{\sqrt{t}} + O(\frac{1}{t})
%\end{align}
%as $t\to\infty$. %Such an expansion is often called an ``asymptotic'' expansion but there is no real need for another name (let $\e=1/\sqrt{t}$ so that $\e\to0$ as $t\to\infty$, which is clearly recognizable as a perturbation expansion).

The of question of which number $N$ of terms is best can be answered easily: $N$ is best when the magnitude of the residual $|\Delta|$ is smallest. Then the error $|z-x|$ satisfies the inequality
\begin{align}
|z-x| \leq \kappa|\Delta|
\end{align}
and $|z-x|$ will also be smallest. Here, $\kappa$ is the condition number of the problem \citep[see, e.g.,][]{CorlessFillion(2013)}, a quantity that is essentially a bound on
\begin{align}
\left|\left(\frac{\partial\Phi}{\partial x}\right)^{-1}\right|\>.
\end{align}
The concept of conditioning of a problem was introduced by Wilkinson \citep[see, e.g.,][]{Wilkinson(1971)}.

%
%
%\subsection{The irrelevance of divergence and the two limits.}
%
%Perturbations seeks to provide valid representations only \emph{locally}.
%\begin{itemize}
%\item Two limits.
%\item Local descriptions and ``patchwork''; matching patches.
%\item The idea of analytic continuation.
%\end{itemize}

%\[
%\lim_{N\to\infty} \sum_{k=0}^N a_k\phi_k(\e) \qquad\qquad \lim_{\e\to0} \sum_{k=0}^N a_k\phi_k(\e)
%\]
%For the second, use a magnifying glass figure as in my slides for 110. This would emphasize locality.
%
%
%\subsection{Formal perturbation series}
%
%Introduce the idea of a formal series --- in particular formal power series, and mention that when the function is sufficiently differentiable this series just happen to be the better known Taylor series.

%\subsection{Patchwork and homotopy}
%
%
%\pagebreak

\section{Singular Perturbations}

A regular perturbation problem may be regarded as one for which the character of the solution set to $\Phi(f,0)=0$ is the same as the character of the solutions to $\Phi(f,\e)=0$ for small $\e$. A singular perturbation problem is one for which the character changes discontinuously at $\e=0$ (\textit{mutatis mutandis}, $t=\infty$). The position of the perturbation parameter may play a key role in determining whether a perturbation is regular or singular, as these examples show:
\[
\begin{array}{c@{\quad}|@{\quad}c}
\textrm{Regular} & \textrm{Singular}\\
\hline
\phantom{\displaystyle \sum_{k=0}^n} x^2-1-\e=0 & \e x^2+x-1=0\\[.25cm]
y''+\e y'+y=0 & \e y''+y'+y=0\\
y(0)=1, y(1)=-1 & y(0)=1,y(1)=-1
\end{array}
\]
There are various definitions of singular perturbation problems in the literature, such as the one by \citet[p.~324]{Bender(1978)}: ``We define a singular perturbation problem as one whose perturbation series either does not take the form of a power series, or, if it does, the power series has a vanishing radius of convergence.'' They also use the more suggestive language:
\begin{quote}
[In either case], the exact solution for $\e=0$ is \emph{fundamentally different in character} from the ``neighbouring'' solutions obtained in the limit as $\e\to0$.
\end{quote}
But this slippery distinction may not always be critical for understanding. Indeed, one common technique for solving so-called singular perturbation problems is to transform them (e.g., by change of variable) into regular perturbation problems. %We will see illustrative examples.
Interestingly, following this characterization, problems with spurious secular terms in the regular expansion display some aspects of both regular and singular perturbation problems (more on this in the example of equation \ref{secterms}).

As a first singular example, consider $\Phi(x,\e) = \e^4 x^5+x-1=0$. This is a prototypical singular perturbation problem. For $\e\neq0$, there are five (complex) roots. However, for $\e=0$, there is only one, namely $x_0=1$. The character---the number of roots---has changed.

We may find the root near $x=1$ by the same method as the previous example: let $z(\e)=a_0+a_1\e+a_2\e^2+a_3\e^3+a_4\e^4+a_5\e^5$ and set as many coefficients of powers of \e\ in $\Delta = \e^4 z^5+z-1$ to zero as possible. This gives $z=1-\e^4+\cdots$ (implicitly stopping at some finite $N$). The other roots are more interesting. To find perturbation series for them, we transform to a regular problem, by using the change of variable  $y=\e x$. %[The $\e^4$ in this example was chosen for simplicity of exposition.]\footnote{But we could use it to explain dominant balance\ldots }
Multiplying the original equation by $\e$ gives
%\[\e(\e^4 x^5 + x-1)=\e\cdot 0\]
%or
\begin{align}
\e^5x^5+\e x-\e=0
\end{align}
and doing the change of variable gives
\begin{align}
y^5 + y - \e = 0\>.
\end{align}
This is a regular perturbation problem for $y$, with five roots. By the same method as before, perturbation series for each can be found:
\begin{align}
y_k(\e) = e^{(2k+1)\pi i/4} + \frac{1}{4e^{(2k+1)\pi i}+1}\e + O(\e^2)
\end{align}
for $0\leq k\leq 3$ and $y_4(\e)=\e-\e^5+O(\e^9)$, giving
\begin{align}
x_k(\e) = \frac{e^{(2k+1)\pi i/4}}{\e} + \frac{1}{4e^{(2k+1)\pi i}+1} + O(\e)
\end{align}
and
\begin{align}
x_4(\e)=1-\e^4+O(\e^8)
\end{align}
as perturbation series solutions for the original (singular) perturbation problem. Notice in this case that four of the solutions go to (complex) infinity as $\e\to0$, leaving only the root near 1 being finite. This observation gives us some geometrical understanding of the reason for which this perturbation problem is singular.

Now, let us examine a problem that turns out to be well-behaved even though we may superficially expect a singular behaviour. Consider
\begin{align}
I(\e) = \int_0^\infty \frac{e^{-x}}{1+\e x}dx\>.
\end{align}
This integral exists for $\e\geq 0$. $I(\e)$ can be expressed in terms of special functions (so-called exponential integrals) but a simpler expression of its behaviour for small $\e$ would be helpful. Naive manipulations give
\begin{align}
\frac{1}{1+\e x} = 1 -\e x + (\e x)^2+\cdots
\end{align}
by the geometric series, so that
\begin{align}
I(\e) &= \int_0^\infty (1-\e x +\e^2x^2+\cdots)e^{-x}dx\>,
\end{align}
suggesting that
\begin{align}
I(\e) & \doteq \int_0^\infty e^{-x}dx -\e\int_0^\infty xe^{-x}dx+\e^2\int_0^\infty x^2e^{-x}dx +\cdots\nonumber\\
&= 1-1!\e +2!\e^2  +\cdots
\end{align}
It may seem remarkable that these naive manipulations, walking over the precipice of divergence for $|\e x|>1$ as they do, manage to give useful perturbative expansions for $I(\e)$. Let us do these manipulations more carefully to see why. Note that
\begin{align}
\frac{1}{1+\e x} = \frac{1-(-\e x)^{N+1}}{1-(-\e x)} + \frac{(-\e x)^{N+1}}{1+\e x}
\end{align}
for all $x$ such that  we avoid the singularity as $\e x= -1$; if we restrict $\e>0$ then $1+\e x\geq 1$ for $x\geq 0$ so division by zero never occurs. Now,
\begin{align}
\frac{1-(-\e x)^{N+1}}{1-(-\e x)} = 1-\e x+\e^2 x^2-\cdots + (-\e x)^N
\end{align}
is the sum of the \emph{finite} geometric series, so
\begin{align}
I(\e) &= \int_0^\infty \left(\sum_{k=0}^N (-1)^k\e^kx^k\right)e^{-x}dx + \int_0^\infty \frac{(-\e x)^{N+1}e^{-x}}{1+\e x} dx\nonumber\\
&= \sum_{k=0}^N (-1)^k\e^k\int_0^\infty x^ke^{-x}dx + \int_0^\infty (-\e)^{N+1}\frac{x^{N+1}e^{-x}}{1+\e x}dx\nonumber\\
&= \sum_{k=0}^N (-1)^k\e^k k! + (-1)^{N+1}\e^{N+1}\int_0^\infty \frac{x^{N+1}e^{-x}}{1+\e x}dx
\end{align}
where all manipulations are of finite quantities and hence legitimate. Now, we do \emph{not} let $N\to\infty$; rather, we look at $\e\to 0$. The term $\int_0^\infty \frac{x^{N+1}e^{-x}}{1+\e x}dx$ is bounded (by $(N+1)!$) and hence the error term is $O(\e^{N+1})$ as $\e\to0$, a term of higher order than all retained terms. Thus this finite sum is a valid perturbation expression for $I(\e)$, in spite of the fact that if we were foolish enough to let $N\to\infty$ we would have a divergent series unless $\e=0$, as can easily be seen for example by the ratio test.

As a practical matter, $I(0.01)$ can be evaluated quite accurately using only three terms of the series:
\begin{align}
I(0.01)\doteq 1-001+2\cdot0.01^2 \doteq 0.9902
\end{align}
The exact value is $I(0.01)=0.9901942287$, so we see that we have obtained a very accurate answer with very little work. This exemplifies why convergence as $N\to\infty$ is of little or no importance in perturbation theory; what matters is that %the perturbation series converges fast in the neighbourhood of the expansion point.
the finite perturbation series be continuous at $\e=0$; indeed differentiable. Of course, this series is useless for large values of $\e$.

%\vspace{1in}
%Perturbation is one of the backbones of applied mathematics. It can be sued to give accurate quantitative predictions from structured nonlinear models, especially near points of interest: at infinity (i.e., for long times), near the starting point, or near to singularities or bifurcation points. The fundamental idea of the method is simple: replace the nonlinear problem $\Phi(f,e)=0$ by a linearized version near $f_0$, i.e., $0=\Phi(f_0,\e)+\Phi'(f_0,\e)\cdot f_1$ and iteratively solve for the successive corrections (this is equivalent to expanding the residual) but the variations seem endless.
%
%It can be used repetitively, in what is known as analytic continuation---this underlies numerical methods. It can be used for infinite diminutional problems. It can be used on the dimension of a problem itself (this is operator regularization). Given a small crack in a problem, perturbation methods can be used to widen that crack, often laying the solution as bare as one could want. We refer the reader to B\&O, O'Malley, van Dyke, etc. for further information.

%\paragraph{Linsdtedt and multiscale}
%
%Probably just cite.
%
%\paragraph{Renormalization}
%
%Probably just cite.
%
%
%\paragraph{Perturbation methods in quantum mechanics}
%
%Just mention its central in QM and that they underly the method of Feynmann diagrams.

%\section{Extra material}

Finally, we discuss an example that illustrates the fact that problems whose perturbative solutions contain secular terms share features of both regular and singular perturbation problems. We borrow this example from \citet[p.~327]{Bender(1978)}. The equation
\begin{align}
y'' + (1-\e x)y = 0 \label{secterms}
\end{align}
with initial conditions $y(0)=1,y'(0)=0$ is an example of a problem that is a \emph{regular} perturbation problem on intervals $0\leq x\leq X$ where $X=O(\e^{-\frac{1}{2}})$, that is, such that $X$ is smaller than constant times $1/\sqrt{\e}$ when $\e$ is small. Indeed, a regular perturbation expansion gives
\begin{multline}
y(x) = \cos(x)+\e\left(\frac{1}{4}x\cos(x) + ( -\frac{1}{4}+\frac{1}{4}x^2) \sin(x) \right) \\
+ \e^2\left(\frac{1}{96}(21x^2-3x^4) \cos(x) +\frac{1}{96}( -21x+10x^3) \sin(x) \right) \\
  +\e^3 \left( (\frac{35}{192}x^3-\frac{35}{128}x-\frac {7}{384}x^5) \cos(x) + (\frac{35}{128}+\frac {7}{96}x^4-\frac{35}{128}x^2-\frac{1}{384}x^6) \sin(x) \right) \\ +O(\e^4)
\end{multline}
However, inspection of this solution reveals that the residual of the truncated series is
\begin{multline}
\Delta = \frac{1}{384}x( -105\sin x -70x^3\cos x-28x^4\sin x +105x\cos x\\
 +105 x^2\sin x+6x^5\cos x+x^6\sin x ) \e^4 +O(\e^5)\>,
\end{multline}
which contains a term $\e^4x^7\cos(x)$ when expanded. Since $z''+(1-\e x)z-\Delta=0$, this is a \emph{small} perturbation of the original equation only when $\e^4x^7\cos(x)$ is small compared to $-\e x\cos(x)$; hence $\e^3x^6\ll 1$ or $x\ll \e^{1/2}$. However, outside of this interval, the effects of the perturbation can be arbitrarily large, and so the problem would be considered singular in that region.

By using various methods that consist in rescaling the problem to eliminate the secular terms, such as Lindstedt's method, renormalization, or multiple scales, we can find an improved solution which is hardly more complicated,
\begin{align}
z = e^{\e x/4}\cos\left(x-\frac{\e t^2}{4}\right) + O(\e^2)\>,
\end{align}
and it is valid for $0\leq x\leq O(1/\e)$. Since $\e^{-1/2}\ll 1/\e$, this is an improvement.
Some applied mathematicians think of these methods as being appropriate for `singular' perturbations, and consider this and problems like it to be `singular'. After all, the residual with the result from the regular perturbation method is not uniformly small as $\e\to0$; it is only small on intervals not as large as $O(\e^{-1/2})$. Notice also that if $|\e x|>1$, then the perturbation is itself not small.

\section{Limitations}

Perturbation methods are extremely powerful and apply in an incredibly varied array of situations. Yet, they have some limitations. A basic limitation of perturbations methods is that they do not apply to discontinuous problems; at least, not unless they have been regularized (transformed into smoother problems). Also, the cost of computing the series may grow quite rapidly with the order of the expansion. We are aware of cases with cost growing exponentially with $N$ if done naively; %this is prohibitive since $O(N^4)$, though polynomial, means in practice that the number of terms that can be computed even on a modern machine may be less than what is wanted.
in many cases this can be improved to polynomial cost, such as $O(N^4)$, but even so this can be expansive.

Perhaps more interestingly, the results are valid only in a (possibly very small) region about the expansion point ($\e=0$); this can give sometimes only a very limited view of the full range of behaviour of the model. But the main limitation is the fact that to find $p_0$ requires solving a nonlinear problem. In general this is always the hardest part. To use perturbation methods, one has to have an initial approximation $p_0\phi_0(\e)$; the problem of finding each successive $p_k$ for $k\geq 1$ thereafter is (usually) linear. These two limitations may be addressed separately using a technique known as homotopy, or continuation. This can be a numerical technique but is not necessarily so, and it is helpful in a broad collection of situations \citep[see, e.g.,][]{Sommese(2005)} %; one can use perturbation about another point, to get started

%[ Put figure here ]

%Consult e.g. Sommese \& Wampler for an introduction. See Rand \& Armbruster for an interesting computer algebra solution of differential equations starting from a problem whose solution is expressed in terms of elliptic functions. Expanding [?]'s library of exact solutions also expands the reach of perturbation.

%\section{Old notes}
%
%If on perturbation:
%
%Cite Abel and Heavyside.
%
%Explain that philosophers of sci have some knowledge of what they are, but it's badly understood why they are such an effective tool.
%
%1. How it works.
%2. Relation to idealization: needs to be kick-started.
%3. Latches onto the world based on a notion of selective accuracy.
%4. Why sometimes we need to take a number of terms that isn't too big (problem of secular terms).
%5. Residual analysis: no need of comparing approximations to experimental results (so, no need to say it's unphysical because unbounded).
%5. Many sets of modelling assumptions can have the same perturbation expansion. (e.g., pendulum and masses on springs are similar because their 1st-order approximation is a simple harmonic oscillator).
%
%
%
%
%
%Explain how their recursive character makes them particularly important for scientific computing and for computer simulations.
%Demonstrate with Euler method and an RK-method.

%\bibliography{../../../bibliography,../../../bibliography2,../../../bibliography3,../../../bibliography4}
%\bibliographystyle{apa}

\end{document}